\documentclass[12pt]{amsart}
\usepackage{amscd,amssymb}
\usepackage[arrow,matrix]{xy}

\theoremstyle{plain}
\newtheorem{thm}[subsection]{Theorem}
\newtheorem{lem}[subsection]{Lemma}
\newtheorem{prop}[subsection]{Proposition}
\newtheorem{cor}[subsection]{Corollary}

\theoremstyle{definition}
\newtheorem{rk}[subsection]{Remark}
\newtheorem{definition}[subsection]{Definition}
\newtheorem{ex}[subsection]{Example}

\numberwithin{equation}{section}
\newcommand{\OO}{{\mathcal O}}
\newcommand{\X}{{\mathcal X}}

\newcommand{\RR}{{\mathcal R}}

\newcommand{\F}{{\mathcal F}}
\newcommand{\A}{{\mathcal A}}
\newcommand{\B}{{\mathcal B}}
\newcommand{\wC}{\widehat{{C}}}

\newcommand{\CC}{{\mathcal C}}
\newcommand{\LL}{{\mathcal L}}

\newcommand{\V}{{\mathcal V}}

\newcommand{\Z}{\mathbb{Z}}
\newcommand{\Q}{\mathbb{Q}}

\newcommand{\C}{\mathbb{C}}

\newcommand{\PP}{\mathbb{P}}
\newcommand{\FF}{\mathbb{F}}
\newcommand{\T}{\mathbb{T}}

\DeclareMathOperator{\Hom}{Hom}

\DeclareMathOperator{\im}{im}

\DeclareMathOperator{\mult}{mult}

\DeclareMathOperator{\codim}{codim}

\begin{document}

\title [Pencils of plane curves and characteristic varieties]
{ Pencils of plane curves and characteristic varieties  }

\author[Alexandru Dimca]{Alexandru Dimca }
\address{  Laboratoire J.A. Dieudonn\'e, UMR du CNRS 6621,
                 Universit\'e de Nice Sophia-Antipolis,
                 Parc Valrose,
                 06108 Nice Cedex 02,
                 FRANCE.}
\email
{dimca@math.unice.fr}



\begin{abstract}
We give a geometric approach to the relation between the irreducible components of the
characteristic varieties
of local systems on a plane curve arrangement complement and the associated 
pencils of plane curves. In the case of line arrangements, this relation  was  recently  discovered   by M. Falk and S. Yuzvinsky \cite{FY} and the geometric point of view was already hinted at by A. Libgober and  S. Yuzvinsky, 
see \cite{LY}, Section 7. Our study yields new geometric insight on the translated components of the
characteristic varieties relating them to the multiplicities of curves in the associated pencil, in a close analogy to the compact situation treated by A. Beauville  \cite{Beau}.


\end{abstract}

\maketitle

\section{Introduction} \label{sec:intro}

Let $\A$ be a line arrangement in the complex projective plane $\PP^2$
and let $M(\A)$ be the corresponding complement. The characteristic varieties 
$\V_m(M(\A))$ (resp. the resonance varieties $\RR_m(M(\A))$) describe the jumping loci
for the dimension of the twisted cohomology group $H^1(M(\A),\LL)$, with $\LL$ a rank one local system on
the complement $M(\A)$, (resp. a rank one local system $\LL$ close to the trivial local system), see for details section 3 below.

 Recently  M. Falk and S. Yuzvinsky \cite{FY} have shown that the existence of a {\it global}
$d$-dimensional irreducible component $E$ in   $\RR_1(M(\A))$ with $d\ge 2$ is equivalent to the existence of a pencil $\CC$ of plane curves on $\PP^2$ with an irreducible generic member such that

\noindent (i) the pencil $\CC$ has $d+1$ fibers $\CC_b$, for $b \in B$ a finite subset of $\PP^1$ with $|B|=d+1$, each one of them
being the union of lines in $\A$ (possibly with some multiplicities);

\noindent (ii) these $d+1$ degenerate fibers $\CC_b$ correspond to a partition of the set of lines
in $\A$.

\noindent  We say that an arrangement of this type is {\it minimal} with respect to the pencil $\CC$
and the set $B \subset \PP^1$, see Definition \ref{def0} below.

\medskip

This surprising equivalence is established in \cite{FY} via a combinatorial approach, based on the description of the irreducible components of the resonant variety  $\RR_1(M(\A))$ in terms of generalized Cartan matrices obtained by Libgober and Yuzvinsky  \cite{LY}.

\medskip

In this paper, the first aim is to reprove this result in the more general setting of curve arrangement complements
$M=\PP^2 \setminus C$, which allows us to grasp the main features of the situation, and, in the end, to better understand even the special case of line arrangements. Indeed, on one hand, a general fiber of a pencil $\CC$  associated to a line arrangement is a curve. On the other hand this setting is more flexible, and we can construct easily vivid examples which are
hard to find in the class of line arrangements, see Example \ref{exfin3} at the end.

The main technical tool is provided by  Arapura's results on the  irreducible components of the characteristic  variety  $\V_1(M)$,  see \cite{A}, section V and Theorem \ref{thm1.5} below. To pass from  the irreducible components of the characteristic variety
to  the irreducible components of the resonance variety, we use one of the main result in
\cite{DPS} (see D. Cohen and A. Suciu \cite{CS1} in the case of line arrangements).
Though our Theorem \ref{thm2} is not as precise as the
description of the line arrangements in  \cite{FY}, we feel that our geometric approach
brings light to what would be otherwise a mysterious property. In fact, we can recover most of the additional results in  \cite{FY} concerning the combinatoric of the arrangement under the additional hypothesis that all the irreducible curves in our arrangement are smooth and intersecting transversally. In the general case, at each base point one has a pencil of plane curve singularities which can be studied, see for instance  \cite{LW} and \cite{DM}.

\medskip

The second and {\it main aim} of this paper is to study the translated components of the
characteristic varieties  $\V_1(M)$. According to  Arapura's results, such a component $W$ is described by
a pair $(f, \rho)$ where

\noindent (a) $f$ is a surjective morphism $M\to S=\PP^1 \setminus B$, which is nothing else but a {\it not necessarily minimal
arrangement} with respect to a given pencil $\CC$ and the set $B$, see Proposition  \ref{prop2}.

\noindent (b)  $\rho$ is a torsion character such that $W$ is the translate by  $\rho$ of a  subtorus constructed via $f$.

Our results can be described briefly as follows:

\noindent (A) The set of mappings $f$ arising in (a) above are parametrized by the (rationally defined) maximal isotropic linear subspaces 
$E \subset H^1(M,\C)$. In fact, in the case $\dim E =1$, not all rationally defined maximal isotropic linear subspaces yield components in $\V_1(M)$, see Remark \ref{rk0.2} (ii).

When $\dim E \ge 2$, then this maximal isotropy condition is equivalent to asking $E$ to be an
irreducible component of
the resonance variety $\RR_1(M)$, see Corollary \ref{cor0.1}, Corollary \ref{cor1.0},
Proposition \ref{prop2.5} and Corollary \ref{cor1.1}. Moreover in this case the rationality condition is automatically fulfilled, see Remark \ref{rk0.1}.

If the arrangement $C$ is given (e.g. the equations $f_j=0$ for the components $C_j$ of $C$ are known), then the map $f$ associated to a (rationally defined) maximal isotropic
linear subspace $E \subset H^1(M,\C)$ can be constructed explicitly, see Propositions \ref{prop4.5} and \ref{prop4.6}.
These results can be regarded as a non-proper Castelnuovo-De Franchis Lemma, see \cite{Beau}, \cite{Cat}.
However, it is not clear whether this construction is combinatorial in the case of a line arrangement.

\noindent (B) The characters $\rho$ arising in (b) above for a given map $f$ are parametrized by the Pontrjagin dual
$T(f)^*=\Hom(T(f), \C^*)$ of a finite group $T(f)$ defined in terms of the topology of the
mapping $f$. This group depends on the multiple fibers in the pencil associated to $f$, see the formulas \eqref{t}, \eqref{k2}, Theorem \ref{thm5} and Corollary \ref{cor5}.
For components of dimension at least 2, any character in $T(f)^*$ actually gives rise to
a component, see Proposition \ref{prop2.5}, while for 1-dimensional components
one should discard the trivial character in $T(f)^*$, see Proposition \ref{prop6} 
which covers a rather general situation.

The 1-dimensional case  is the most mysterious, and Suciu's example of such a component
for the deleted $B_3$-arrangement given in  \cite{S1}, \cite{S2}
played a key role in our understanding of this question. We consider this component in detail in Example \ref{exMASTER} and at the end of the paper, together with its generalization given by the $\A_m$-arrangements discussed in \cite{C1} and \cite{CDS},  as a good test for our results.

\medskip
Completely similar results on the characteristic varieties of rank 1 local systems on a compact K\"ahler manifold
were obtained by A. Beauville in \cite{Beau}. The techniques of proof are rather different and it does not seem easy to obtain our results by using Beauville's.

\medskip
In section 2 we collect some basic facts on rational mappings  $f:\PP^2 \to \PP^1$
and the associated pencils. Lemma \ref{lem1} intends to clarify the key notion of admissible map used by Arapura in  \cite{A}.

In section 3 we give the main definitions and properties of characteristic and resonant varieties. 
Theorem  \ref{thm1.5} collects some (more or less known) facts on the irreducible components of the characteristic varieties, which are derived by a careful reading of Arapura's paper  \cite{A}. We also prove the major part of the claims in (A).

 In section 4,  Theorem \ref{thm2} is our analog of the main results in \cite{FY}. As a consequence, we obtain a necessary numerical condition involving self-intersection numbers for the existence of an essential  positive dimensional irreducible component of a characteristic variety. This is practically the same condition as that in Theorem 4.1.1 in Libgober  \cite{L1}, which was established via the use of adjunction ideals.

In section 5 we discuss the complements $M$ which are in an obvious way  $K(\pi, 1)$-spaces,
i.e. for which the mapping $f:M \to S$ considered above is a fibration, and we conclude with Example \ref{ex2}
where several of the above features are clearly illustrated. In particular we point out several differences with the case of  line arrangements.

In the final section we associate to a map  $f:M \to S$ as above a finite abelian group $T(f)$, such that the 
 torsion character $\rho$ is determined by a character $\tilde \rho$ of $T(f)$, see formula  \eqref{k2}.
Then we compute this  group $T(f)$ in terms of the multiplicities of some special fibers of the pencil associated to $f$, see Theorem \ref{thm5}. This result explains why usually $T(f)=1$ and hence there are no translated components
associated to $f$. Several key examples are also included here.

A list of very interesting related open questions may be found in Libgober
 \cite{L2}.

We would like to thank Alexander Suciu for interesting and stimulating discussions on the subject of this paper, and for suggesting several improvements
of the presentation.

\section{On rational maps from $\PP^2$ to $\PP^1$} \label{sec:two}

Let $f:\PP^2 \to \PP^1$ be a rational map. Then there is a minimal non-empty finite set $A \subset \PP^2$
such that $f$ is defined on $U=\PP^2 \setminus A$. We recall the following basic fact. 

\begin{prop} \label{prop1}
Any  morphism $f:U \to \PP^1$ is given by a pencil $\CC: \alpha_1 P_1+\alpha _2 P_2$ of plane curves
having the base locus $V(P_1,P_2)$ the minimal finite set $A$. This pencil is unique up-to an isomorphism of
$\PP^1$. 
\end{prop} 

\proof  We have to show the existence of two homogeneous polynomials $P_1,P_2 \in \C[X,Y,Z]$ of same degree
$D$, called the degree of the pencil, such that 

\begin{enumerate}

\item $V(P_1,P_2)=\{x \in \PP^2~|~~P_1(x)=P_2(x)=0\}= A$ (in particular, these polynomials have no common factor), and

\item for any point $x \in U$, one has $f(x)=(P_1(x):P_2(x))$.

\end{enumerate}

It is well known, see for instance \cite{H}, p. 150, that a morphism $f:U \to \PP^1$ is given
by a line bundle $\LL \in Pic(U)$ and two sections $s_1, s_2 \in \Gamma (U,\LL)$ which do not vanish
both at any point in $U$. In fact $\LL=f^*(\OO(1))$ and $s_i=f^*(y_i)$, with $y_1,y_2$ a system of homogeneous coordinates on $\PP^1$. With this notation, one has $f(x)=[a:b]$ where $[a:b] \in \PP^1$ is such that
$as_2(x)-bs_1(x)=0$.

Since $U$ is smooth, we have $Pic (U) = C\ell (U)$ and similarly
$Pic ( \PP^2 ) = C\ell ( \PP^2)$, see for instance \cite{H}, p. 145. On the other hand, the inclusion
$j:U \to \PP^2$ induces an isomorphism $j^*:C\ell(\PP^2) \to C\ell (U)$, as $\codim A =2$, see \cite{H}, p. 133. It follows that $j^*:Pic(\PP^2) \to Pic (U)$ is also an isomorphism, i.e. any line bundle $\LL \in Pic(U)$
is the restriction to $U$ of a line bundle $\OO(D)$ and the global sections of $\LL$ are nothing else but the restrictions of global sections of the line bundle $\OO(D)$, which are the degree $D$ homogeneous polynomials.

\endproof

Let $C \subset \PP^2$ be a reduced curve such that $C=\cup _{j=1,r}C_j$, with $C_j$ irreducible
curve of degree $d_j$. We set $M=\PP^2 \setminus C$.

\begin{prop} \label{prop2}
Let $B \subset \PP^1$ be a finite set and denote by $S$ the complement $\PP^1\setminus B$.
For any surjective morphism $f:M \to S$ there is a pencil $\CC: \alpha_1 P_1+\alpha _2 P_2$ of curves
in $\PP^2$ such that any irreducible component $C_j$ of $C$  is in one of the three following cases.

\begin{enumerate}

\item $C_j$ is contained in a curve $\CC_b$ in the pencil $\CC$, corresponding to a point $b \in B$;

\item $C_j$ is strictly contained in a curve $\CC_s$ in the pencil $\CC$, corresponding to a point $ s \in S$;

\item $C_j$ is a horizontal component, i.e. $C_j$ intersects the generic fiber $\CC_t$ of the pencil $\CC$ outside the base locus.

\end{enumerate}

Moreover, $C_j$ is in the first case above if and only if the homology class $\gamma_j$ of a small loop around $C_j$
satisfies $H_1(f)(\gamma_j) \ne 0$ in $H_1(S,\Z)$.

\end{prop} 

\proof Let $i:S \to\PP^1$ be the inclusion and set $g=i\circ f: M\to \PP^1$. Then $g$ is a rational map 
 $\PP^2 \to \PP^1$ and there is a minimal non-empty finite set $A \subset \PP^2$
such that $g$ can be extended to a morphism $g_1:U=\PP^2 \setminus A \to \PP^1$. We apply to $g_1$ 
Proposition \ref{prop1} and get the corresponding pencil $\CC$. Let $C_j$ be an irreducible component of $C$.
Then either $g_1(C_j)$ is a point, which leads to the first two cases, or $g_1(C_j)$ is dense in $\PP^1$, which leads to the last case. The strict inclusion in the second case comes from the surjectivity of $f$.

The last claim is obvious. For instance, in the first case, if $\delta _b$ is a small loop at $b$, then one has $H_1(f)(\gamma_j)=m_j\cdot \delta _b$, with $m_j>0$ the multiplicity of the curve
$C_j$ in the divisor $g_1^{-1}(b)$ (if the orientations of the loops $\gamma_j$ and $\delta _b$ are properly chosen).

\endproof

\begin{definition} \label{def0}
 In the setting of Proposition \ref{prop2}, we say that the curve arrangement $C$ is {\it minimal} with respect to the surjective mapping $f:M\to S$ if any component $C_j$ of $C$ is of type (1), i.e.
 $C_j$ is contained in a curve $\CC_b$ in the pencil $\CC$, corresponding to a point $b \in B$.
We say that the curve arrangement $C$ is  {\it special} with respect to the surjective mapping $f:M\to S$ if some component $C_j$ of $C$ is of type (2), i.e. $C_j$ is strictly contained in a curve $\CC_s$ in the pencil $\CC$, corresponding to a point $ s \in S$.

\end{definition} 

\begin{rk} \label{rk1}

If $|B|>1$, then the base locus $\X$ of the pencil $\CC$ is just the intersection of any two distinct fibers
$\CC_b \cap \CC_{b'}$ for $b,b'$ distinct points in $B$. Note also that the second case cannot occur if all
 the fibers $\CC_s$ for $s \in S$ are irreducible. Note that  the fibers $\CC_s$ may be non-reduced, i.e. we consider them usually as divisors. Saying that $C_j$ is contained in $\CC_s$ means that $\CC_s=m_j C_j+...$,
with $m_j>0$. On the other hand, $C$ is a reduced curve.

\end{rk} 

We conclude this section by the following easy fact.

\begin{lem} \label{lem1}
Let $X$ and $S$ be  smooth irreducible algebraic varieties, $\dim S=1$ and let $f:X \to S$ be a non-constant
morphism. Then for any compactification $f':X' \to S'$ of $f$ with $X'$, $S'$ smooth, the following are equivalent.

\medskip

\noindent (i) The generic fiber $F$ of $f$ is connected.

\medskip

\noindent (ii) The generic fiber $F'$ of $f'$ is connected.

\medskip

\noindent (iii) All the fibers  of $f'$ are connected.

\medskip

If these equivalent conditions hold, then $f_{\sharp}: \pi _1(X) \to \pi_1(S)$ and $f_{\sharp}': \pi _1(X') \to \pi_1(S')$ are surjective.

\end{lem}

\proof Note that $D=X' \setminus X$ is a proper subvariety (not necessarily a normal crossing divisor) with finitely many irreducible components $D_m$. For each such component $D_m$, either $f'(D_m)$ is a point,
or $f':D_m \to S'$ is surjective. In this latter case, it follows that $\dim (F' \cap D_m) < \dim D_m \leq \dim F'$. Since $F'$ is smooth of pure dimension, it follows that $F'$ is connected if and only if 
$F =F' \setminus \cup_m ( D_m \cap F')$ is connected.
 To show that $(ii)$ implies $(iii)$ it is enough to use the Stein factorization theorem,
  see for instance \cite{H}, p. 280, and the fact that a morphism between two smooth projective curves which is of degree one
 (i.e. generically injective) is in fact an isomorphism.

\medskip
 
 To prove the last claim for $f$, note that there is a Zariski open and dense subset $S_0 \subset S$ such that
 $f$ induces a locally trivial topological fibration $f:X_0=f^{-1}(S_0) \to S_0$ with fiber type $F$. Since $F$ is connected, we get an epimorphism $f_{\sharp}: \pi _1(X_0) \to \pi_1(S_0)$. The inclusion of $S_0$ into $S$ induces an epimorphism as well at the level of fundamental groups. Let $j:X_0 \to X$ be the inclusion. Then
 we have seen that $f \circ j$ induces an epimorphism as well at the level of fundamental groups. Therefore the same is true for $f$. The proof for $f'$ is completely similar.

 \endproof
 
 \section{Local systems, characteristic varieties, and resonance varieties} \label{sec:3}

 \subsection{Local systems on $S$} \label{sec:3.1}
 
 Here we return to the notation $S=\PP^1\setminus B$, with $B=\{b_1,...,b_k\}$ a finite set of cardinal
 $|B|=k>0$. If $\delta_i$ denotes an elementary loop based at some base point $b \in B$ and turning once around the point $b_i$, then using the usual choices, the fundamental group of $S$ is given by
 
\begin{equation} \label{eq1} 
\pi_1(S)=<\delta _1,...,\delta_k ~ | ~ \delta _1\cdots \delta_k=1>.
\end{equation}

It follows that the first integral homology group is given by
\begin{equation} \label{eq2} 
H_1(S)=\Z <\delta _1,...,\delta_k>/<\delta _1+...+ \delta_k=0>.
\end{equation}
Therefore, the rank one local systems on $S$ are parametrized by the $(k-1)$-dimensional algebraic torus

\begin{equation} \label{eq3} 
\T(S)= \Hom( H_1(S),\C^*)=\{ \lambda=(\lambda_1,...,\lambda_k)\in (\C^*)^k~|~ \lambda _1\cdots \lambda_k=1\}.
\end{equation}
Here $\lambda_j\in \C^*$ is the monodromy about the point $b_j \in B$. For $ \lambda \in \T(S)$, we denote by $\LL_{ \lambda}$ the corresponding rank one local system on $S$.

The twisted cohomology groups $H^m(S, \LL_{ \lambda})$ are easy to compute. There are two cases.

\noindent Case 1 ($\LL_{ \lambda}=\C_S$). Then we get the usual cohomology groups of $S$, namely
we have $\dim H^0(S, \LL_{ \lambda})=1$, $\dim H^1(S, \LL_{ \lambda})=k-1$ and $H^m(S, \LL_{ \lambda})=0$
for $m \geq 2$.

\noindent Case 2 ($\LL_{ \lambda}$ is nontrivial). This case corresponds to the case when at least one monodromy $\lambda_j$ is not 1. Then
we have $\dim H^0(S, \LL_{ \lambda})=0$, $\dim H^1(S, \LL_{ \lambda})=k-2$ and $H^m(S, \LL_{ \lambda})=0$
for $m \geq 2$.

\subsection{Local systems on $M$} \label{sec:3.2}

Let $\gamma_j$ be an elementary loop around the irreducible component $C_j$, for $j=1,...,r.$
Then it is known, see for instance \cite{D1}, p. 102, that
\begin{equation} \label{eq4} 
H_1(M)=\Z <\gamma _1,...,\gamma_r>/<d_1\gamma _1+...+ d_r\gamma_r=0>
\end{equation}
where $d_j$ is the degree of the component $C_j$.
It follows that the rank one local systems on $M$ are parametrized by the algebraic group
\begin{equation} \label{eq5} 
\T(M)= \Hom( H_1(S),\C^*)=\{ \rho=(\rho_1,...,\rho_r)\in (\C^*)^r ~|~  \rho _1^{d_1} \cdots \rho_r^{d_r}=1\}.
\end{equation}
The connected component $\T^0(M)$ of the unit element $1 \in \T(M)$ is the $(r-1)$-dimensional
torus given by
\begin{equation} \label{eq6} 
\T^0(M)=\{ \rho=(\rho_1,...,\rho_r)\in (\C^*)^r ~|~  \rho _1^{e_1} \cdots \rho_r^{e_r}=1\}
\end{equation}
with $D=G.C.D.(d_1,...,d_r)$ and $e_j=d_j/D$ for $j=1,...,r.$

\begin{rk} \label{rk0}
If $d_1=1$, then $\{\gamma_2,..., \gamma_r\}$ is a basis for $H_1(M)$ and the torus $\T(M)$ can be identified to
$(\C^*)^{r-1}$ under the projection $\rho \mapsto (\rho_2,...,\rho_r)$.

\end{rk}

Now the computation of the twisted cohomology groups $H^m(M, \LL_{ \rho})$ is one of the major problems.
The case when $\LL_{ \rho}=\C_M$ is easy, and the result depends on the local singularities of the plane curve
$C$. In fact $\dim H^0(M, \C)=1$, $\dim H^1(M, \C)=r-1$ and $H^m(M, \C)=0$
for $m \geq 3$. To determine the remaining Betti number $b_2(M)=\dim H^2(M, \C)$ is the same as determining
the Euler characteristic $\chi (M)=3-\chi (C)$ and this can be done, e.g. by using the
formula for $\chi(C)$ given  in  \cite{D1}, p. 162. 

In the sequel we concentrate on the case 
$\LL_{ \rho} \ne \C_M$ and assume $\chi (M)$ known. Then we have $H^m(M, \LL_{ \rho})=0$
for $m=0$ and $m \geq 2$, and $\dim H^2(M, \LL_{ \rho})-\dim H^1(M, \LL_{ \rho})=\chi(M)$, see for instance
\cite{D2}, p. 49. To study these cohomology groups, one idea is to study the characteristic varieties
\begin{equation} \label{eq7} 
\V_m(M)=\{ \rho=(\rho_1,...,\rho_r)\in (\C^*)^r ~|~  \dim H^1(M, \LL_{ \rho}) \ge m \}.
\end{equation}

\begin{definition} \label{def1} 

An irreducible component $W$ of such an $m$-th characteristic variety $\V_m(M)$ is called a
 {\it coordinate component} (resp. a  {\it translated coordinate component})
if $W$ is contained in a subgroup $\T_j$ of $\T(M)$ defined by an equality $\rho_j=1$ for some $j$
(resp. there is a torsion character $\rho \in \T(M)$ such that $W \subset \rho \T_j$ for some $j$).
An irreducible component $W$ which is not a translated coordinate component is called a
 {\it global component}.
\end{definition}

Note that if $1 \in W$, then $W$ is a  coordinate component if and only if  $W$ is a  translated coordinate component.

\medskip

 Let $C(j)$ be the
plane curve obtained from $C$ by discarding the $j$-th component $C_j$. Let $M(j)=\PP^2 \setminus C(j)$ be the corresponding complement. Then the inclusion $ \iota_j: M \to M(j)$ induces an epimorphism $H_1(M) \to H_1( M(j))$
and hence an embedding $\iota ^*_j : \T(M(j)) \to \T(M)$. 

\begin{definition} \label{def2} 

An irreducible component $W$ of the $m$-th characteristic variety
$\V_m(M)$ is called a {\it non-essential component, or a pull-back component} if $W =\iota^*_j(W_j)$ for some $j$ and some irreducible component $W_j$ of the $m$-th characteristic variety
$\V_m(M(j))$, see \cite{ACC}, \cite{FY}, \cite{L1}. An irreducible component $W$ which is not  non-essential  is called an {\it  essential component}.
\end{definition}

\begin{rk} \label{rk2}
The notions of coordinate and (non-)essential  component depend on the curve arrangement $C=\cup C_i$, i.e. on the chosen embedding of $M$ into $\PP^2$. So they are not invariants of the surface $M$.
For more on this see  \cite{ACC}.

\end{rk}

Assume given a surjective morphism $f:M \to S$ such that $f_{\sharp}: \pi _1(M) \to \pi_1(S)$ is surjective.
This gives rise to an embedding $f^*:\T(S) \to \T^0(M)$, which implies in particular $k \le r$.
More precisely, if we start with $\LL_{ \lambda} \in \T(S)$, then the monodromy $\rho _j$ of the pull-back
local system $f^*\LL_{ \lambda} =\LL_{ \rho}$ is given by

\medskip

\noindent (i) $\rho _j=1$  if the component $C_j$ is not in the first case of Proposition \ref{prop2}, and by

\medskip

\noindent (ii) $\rho _j=\lambda _i^{m_j}$  if the component $C_j$ is in the first case of Proposition \ref{prop2}, i.e. $g_1(C_j)=b_i$ in the notation from the proof of Proposition \ref{prop2}. Recall that $m_j$ is the multiplicity of $C_j$ in $f^{-1}(b_i)$.

\begin{cor} \label{cor1} 
With the above notation, the pull-back local system $f^*\LL_{ \lambda} =\LL_{ \rho}$
satisfies $\rho_j \ne 1$ for all $j=1,...,r$ if and only if

\noindent (i) The curve $C$ consists exactly of the fibers of the associated pencil $\CC$ 
corresponding to the points in $B$.

\noindent (ii) For all $j=1,...,r$, if we set $g_1(C_j)=b_{i(j)}$, then $\lambda ^{m_j}_{i(j)}\ne 1$.

\end{cor} 

\subsection{Arapura's results} \label{sec:3.3}

We recall here some of the main results from \cite{A}, applied to the rank one local systems on $M$, with some additions from  \cite{L1}, \cite{DPS} and some new consequences.

\begin{thm} \label{thm1.5}

 Let $W$ be an irreducible component of $\V_1(M)$ and assume that  $d_W :=\dim W \ge 1$.
 Then there is a surjective morphism $f_W:M \to S_W$ with connected generic fiber $F(f_W)$,
and a torsion character $\rho_W \in \T(M)$  such that  
$$W= \rho_W \otimes  f_W^*(\T(S_W)).$$
More precisely, the following hold.

\medskip

\noindent (i) $S_W=\PP^1 \setminus B_W$, with  $B_W$ a  finite set  satisfying $k_W:=|B_W|=d_W+1$.

\medskip

\noindent (ii)  For any local system $\LL \in W$, the restriction $\LL | F(f_W)$ of $\LL$ to the generic fiber of $f_W$ is trivial, i.e. $\LL | F(f_W)=\C_{ F(f_W)}$.


\medskip

\noindent (iii) If $N_W$ is the order of the character $ \rho_W$, then there is a commutative diagram
$$\xymatrix{
M'  \ar[r]^p \ar[d]^{f'_W} & M \ar[d]^{f_W} \\
S'_W  \ar[r]^q & S_W
}$$
where $p$ is a unramified $N_W$-cyclic Galois covering , $q$ is possibly ramified  $N_W$-cyclic Galois covering, $f'_W$ is $\mu _{N_W}$-equivariant in the obvious sense, has a generic fiber $F(f'_W)$ isomorphic to the  generic fiber $F(f_W)$ of  $f_W$,
and $p^*\rho_W $ is trivial. Here $\mu _{N_W}$ denotes the cyclic group of the $N_W$-th roots of unity.

\medskip

\noindent (iv) If $1 \in W$ and $\LL  \in W$, then $\dim H^1(M,\LL) \ge -\chi(S)=d_W-1$ and equality holds with finitely many exceptions.

\medskip

\noindent (v) If $1 \in W$, then $d_W \ge 2.$

\medskip

\noindent (vi) If $1 \notin W$ and $d_W \ge 2$, then the subtorus $W'= f_W^*(\T(S_W))$ is another irreducible component of $\V_1(M)$. In this situation, $W'$ is a coordinate component if and only if $W$
is a translated coordinate component.

\end{thm}

\proof

The first claim is just Thm. 1.6 in \cite{A}, section V.

\medskip

Now we prove the claims (i) and (vi). The fact that $S_W$ has to be rational in this situation is noted in \cite{L1}, and it follows
 from the fact that a compactification $\overline M$ of $M$ is simply-connected and hence, via
Lemma \ref{lem1}, we get that  a compactification $\overline S_W$ of $S_W$ is simply-connected as well.
A different proof follows from Proposition 5.10 (2) in  \cite{DPS}. When $1 \in W$, the equality 
$k_W=d_W+1$ was  noted in \cite{L1}, see also  Proposition 6.3 in  \cite{DPS}.

Consider now the situation $1\notin W$. Then there are two cases.

\noindent {\bf Case 1.} $k_W=2$. Then $S_W=\C^*$ and hence $1 \le d_W \le \dim \T(S_W)=1$.

\medskip

\noindent {\bf Case 2.} $k_W \ge 3$. Then $\chi(S_W)<0$ and  $W'= f_W^*(\T(S_W))$ is another irreducible component of $\V_1(M)$ by Prop. 1.7 in \cite{A}. Since $1 \in W'$, $d_{W'}=d_W$, we get the equality
$k_W=d_W+1$. Moreover, this implies $d_W\ge 2$, i.e. we get the claim (vi) as well.

\medskip

Now we prove the claim (ii). Since  $W= \rho_W \otimes  f_W^*(\T(S_W))$, it is enough to prove that 
 $ \rho_W  | F(f_W)=\C_{ F(f_W)}$. And this is proved in the final part of the proof of Prop. 1.3 in 
 \cite{A}. Just note that on the last line of this proof, one should replace ``which forces $\psi |F$ to be trivial''
by  ``which forces $\psi | (F\cap X)$ to be trivial''. (This is due to the fact that $F$ in  \cite{A}
denotes the compactification of our $F= F(f_W)$, and $X$ in  \cite{A} corresponds to our $M$.)

\medskip

The claim (iii) is just the ``untwisting'' part of the proof of Thm. 1.6  in  \cite{A}. The existence of the diagram is explained there via the Stein factorization for $f_W \circ p$. However, the fact that the morphism $q$ has degree $N_W$ depends on the previous claim (ii), and this key point is not mentioned in  \cite{A}.

\medskip

The proof of the claim (iv) is more technical. Using the Projection Formula
\begin{equation} \label{eqp1} 
p_*(\C_{M'})\otimes \LL \simeq p_*(p^*(\LL))
\end{equation}
for $\LL \in W$, see for instance  \cite{D2}, p.42 and then the Leray Spectral Sequence for $p$,
see for instance  \cite{D2}, p. 33, one gets an isomorphism of $\mu _{N_W}$-representations
\begin{equation} \label{eqp2} 
H^1(M',p^*\LL)=H^1(M,p_*(\C_{M'})\otimes \LL ).
\end{equation}
Following the argument in the proof of Thm. 1.6  in  \cite{A}, we get the following
$$\dim H^1(M,\LL) \ge -\chi (S_W)=k_W-2=d_W-1.$$
The only point which deserves some attention is the fact that $S_W$ and $S'_W$ do not admit finite triangulations as claimed  in  \cite{A}, since they are not compact. However, we can replace them by finite
simplicial complexes without changing the topology, e.g. $S_W$ can be replaced by the compact Riemann surface with boundary obtained from $\PP^1$ by deleting small open discs centered at the points in $B_W$.

The fact that there are only finitely many local systems $\LL \in W$ such that $\dim H^1(M,\LL)\ge d_W$ follows
by an argument similar to the end of the proof of Prop. 1.7 in  \cite{A}, section V, see for details
Remark \ref{rk0.15} below.

 \medskip

Finally, the claim (v) is well-known, see for instance  \cite{LY} or  Cor. 6.4 in  \cite{DPS}.

\endproof

Note that the claim (vi) above is obviously false for $d_W=1$ by (v). 
A deeper fact is that (iv) is false when $1\notin W$, see Corollary \ref{cor6}.

\begin{rk} \label{Rk}
Conversely, if $f:M \to S$ is a morphism with a generic connected fiber and with $\chi(S)<0$, then
$W_f=f^*(\T(S))$ is an irreducible component in $\V_1(M)$ such that $1 \in W_f$ and $\dim W_f \ge 2$, see
 \cite{A}, Section V, Prop. 1.7. Some basic situations of this general construction are the following.
 
 \medskip
 
\noindent (i) {\it The local components}, see for instance  \cite{S1}, subsection (2.3) 
in the case of line arrangements. The case of curve arrangements runs as follows. Let $p \in \PP^2$ be a point such that there is a degree $d_p$ and an integer $k_p >2$ such that

\begin{enumerate}

\item the set $A_p=\{j ~~|~~p \in C_j \text{ and } \deg C_j=d_p\}$ has cardinality $k_p$;

\item $\dim <f_j~~|~~j \in A_p>=2,$ with $f_j=0$ being an equation for $C_j$.

\end{enumerate}
If $\{P,Q\}$ is a basis of this 2-dimensional vector space, then the associated pencil induces a map
$$f_p:M \to S_p$$
where $S_p$ is obtained from $\PP^1$ by deleting the $k_p$ points corresponding to the curves
$C_j$, for $j \in A_p$. In this way, the point $p$ produces an irreducible component in $\V_1(M)$, namely
$$W_p= f^*_p(\T(S_p))$$
of dimension $k_p-1$, and which is called local because it depends only on the chosen point $p$. Note that in the case of line arrangements $p$ can be chosen to be any point of multiplicity at least $3$.

\medskip

\noindent (ii) {\it The  components associated to neighborly partitions}, see  \cite{LY}, corresponds
exactly to pencils associated to the line arrangement, as remarked in  \cite{FY}, see the proof of Theorem 2.4
\end{rk}

All these points are illustrated by the following.

\begin{ex} \label{exMASTER} This is a key example discovered by A. Suciu, see Example 4.1 in  \cite{S1}
and Example 10.6 in  \cite{S2}. Consider the line arrangement in $\PP^2$ given by the equation
$$xyz(x-y)(x-z)(y-z)(x-y-z)(x-y+z)=0.$$

We number the lines of the associated affine arrangement in $\C^2$ (obtained by setting $z=1$) as follows: $L_1: x=0$, $L_2:x-1=0$, $L_3: y=0$, $L_4: y-1=0$, $L_5:x-y-1=0$, $L_6:x-y=0$ and $L_7:x-y+1=0$, see the pictures in  Example 4.1 in  \cite{S1} and Example 10.6 in  \cite{S2}.
As stated in  Example 4.1 in  \cite{S1}, there are

\medskip

\noindent (i) Seven local components: six of dimension 2, corresponding to the triple points, and one of dimension 3, for the quadruple point.

\medskip

\noindent (ii) Five components of dimension 2, passing through 1, coming from the following neighborly
partitions (of braid subarrangements): $(15|26|38)$, $(28|36|45)$, $(14|23|68)$, $(16|27|48)$ and $(18|37|46)$. For instance, the pencil corresponding to the first partition is given by
$P=L_1L_5=x(x-y-z)$ and $Q=L_2L_6=(x-z)(x-y)$. Note that $L_3L_6=yz=Q-P$, a fiber in this pencil.

\medskip

\noindent (iii) Finally, there is a 1-dimensional component $W$ in $\V_1(M)$ with 
$$\rho_W=(1,-1,-1,1,1,-1,1,-1) \in \T(M) \subset (\C^*)^8$$
and $f_W:M \to \C^*$ given by
$$f_W(x:y:z)=\frac{x(y-z)(x-y-z)^2}{(x-z)y(x-y+z)^2}$$
or, in affine coordinates
$$f_W(x,y)=\frac{x(y-1)(x-y-1)^2}{(x-1)y(x-y+1)^2}.$$
Then $W \subset \V_1(M)$ and $W \cap  \V_2(M)$ consists of two characters, $\rho_W$ above and
$$\rho _W'=(-1,1,1,-1,1,-1,1,-1).$$
Note that this component $W$ is a translated coordinate component. This is related to the fact that the associated pencil is special. More on this aspect at the end of the paper.

\end{ex}

\medskip

It is clear that any non-essential  component is a  coordinate component.
The following converse result on  positive dimensional coordinate components $W$ of $\V_m(M)$ was obtained by Libgober \cite{L1}. For reader's convenience we include a proof.

\begin{prop} \label{prop3}
Any positive dimensional translated coordinate component of $\V_m(M)$ is non-essential.
\end{prop} 

\proof 
 
 Let $W=\rho_W \otimes f_W^*(\T(S_W))$ be a positive dimensional irreducible component of $\V_m(M)$.  Assume $W$ is contained in the subtorus of $\T^r$ given by $\rho_j=1$. It follows that the corresponding component
 $\rho_{W,j}$ of the character $\rho$ is 1, and that the torus $\T_W= f_W^*(\T(S_W))$ is also contained
 in the same subtorus. 
 The discussion before Corollary \ref{cor1} implies that the corresponding component $C_j$ of $C$ is not in the first case of Proposition \ref{prop2}.
This in turn implies the existence of an extension $f(j): M(j) \to S_W$,
whose generic fibers are still connected (being obtained from those of $f$ by adding at most finitely many points). It follows that $W=\iota^* _j(W_j)$, with $\iota_j:M \to M(j)$ the inclusion and
$W_j=\rho_j \otimes f(j)^*(\T(S_W))$, where the character $\rho_j$ is obtained from  $\rho_W$ by discarding the $j$-th component. To show that $W_j$ is
an irreducible component in $\V_m(M(j))$, we can use Proposition \ref{prop2.5} in the case
$\chi(S_W)<0$ and Corollary  \ref{cor1.1} in the case $S_W=\C^*$ (the fact that $K_{f(j)}=(\iota _j^*)^{-1}(K_f)$ is a rationally defined maximal isotropic subspace is obvious). In the both cases we have to use in addition the equality
$$\dim H^1(M,\LL_1 \otimes f_W^*\LL_2)=\dim H^1(S_W,R^0f_{W*}\LL_1 \otimes \LL_2)=$$
$$=\dim H^1(S_W,R^0f(j)_{*}\LL'_1 \otimes \LL_2)=\dim H^1(M(j),\LL'_1 \otimes f(j)^*\LL_2).$$
Here $\LL_1 \in \T(M)$ (resp. $\LL'_1 \in \T(M(j))$ is the local system corresponding to $\rho_W$
(resp. $\rho_j $), the first and the third equalities come from Remark \ref{rk0.15}, while the middle equality comes from $\LL'_1=\iota_{j,*} \LL_1$ and $f(j) \circ \iota_j=f_W$.

\endproof

In view of this result, it is natural to study first the non-coordinate positive dimensional components.
Indeed, the other components come from simpler arrangements, involving fewer components $C_j$'s. The situation of translated components is different, e.g. the component $W$ studied in Example  \ref{exMASTER} is NOT coming from a simpler arrangement.
The case of 0-dimensional components is very interesting as well, see 
 \cite{ACC},  \cite{S1}, 
 \cite{S2}.

\subsection{Resonance varieties} \label{sec:3.4}

Let $H^*(M,\C)$ be the cohomology algebra of the surface $M$ with $\C$-coefficients. Right multiplication by an element $z \in H^1(M,\C)$ yields a cochain complex $(H^*(M,\C), \mu_z)$. The {\it resonance varieties} of $M$ are the jumping loci for the degree one cohomology of this complex:
\begin{equation} \label{rv1} 
\RR_m(M)=\{z \in H^1(M,\C) ~|~ \dim H^1(H^*(M,\C), \mu_z) \ge m \}.
\end{equation}
One of the main results in  \cite{DPS} gives the following. For the case of hyperplane arrangements see \cite{CS1}.

\begin{thm} \label{thm1}
The exponential map $\exp: H^1(M,\C) \to \T^0(M)$ induces for any $m \ge 1$ an isomorphism of analytic germs
$$(\RR_m(M),0) \simeq (\V_m(M),1).$$
\end{thm}
 The following easy consequence will play a key role.
\begin{cor} \label{cor0.1}
The irreducible components of $\RR_1(M)$ are precisely the maximal linear subspaces $E \subset H^1(M,\C)$, isotropic with respect to the cup product on $M$ 
$$\cup: H^1(M,\C) \times H^1(M,\C) \to H^2(M,\C)$$
 and such that $\dim E \ge 2$.
\end{cor}
\proof

Let $E$ be a component of $\RR_1(M)$.
By the above Theorem there is a component $W$ in $\V_1(M)$ such that $1 \in W$ and $T_1W=E$. 
By Theorem \ref{thm1.5} we can write $W=f^*(\T(S))$, and hence $T_1W=f^*(H^1(S,\C))$ is isotropic
 with respect to the cup product, since the cup product on $H^1(S,\C)$ is trivial.
Maximality of $E$ comes from the fact that $E$ is a component of $\RR_1(M)$. The restriction
$\dim E \ge 2$ comes from  Theorem \ref{thm1.5}, (v).

\endproof

\begin{rk} \label{rk0.1}
It follows from the proof of Corollary \ref{cor0.1} that any  maximal isotropic linear subspaces 
$E \subset H^1(M,\C)$ is {\it rationally defined}, i.e. there is a linear subspace $E_{\Q}\subset H^1(M,\Q)$
such that $E =E_{\Q} \otimes _{\Q}\C$ under the identification $ H^1(M,\C)= H^1(M,\Q) \otimes _{\Q}\C$.
Indeed, one can take $E_{\Q}=f^*(H^1(S,\Q))$.

\end{rk}

We can restate the above Corollary as follows.

\begin{cor} \label{cor1.0}
If $f:M \to S$ is a surjective morphism with connected generic fiber $F$ and $\rho \in \T(M)$ is a torsion character such that  $W=\rho \otimes f^*(\T(S))$ is an irreducible component of $\V_1(M)$ with $\dim W \ge 2$, then
$K_f=f^*(H^1(S,\C))$ is a (rationally defined) maximal isotropic subspace in $H^1(M,\C)$ with respect to the cup-product.
\end{cor}

\proof
Using Theorem \ref{thm1.5}, (vi), we can take $\rho =1$. Then $K_f$ is exactly the tangent space at $1 \in W$,
and, by Theorem \ref{thm1}, an irreducible component of $\RR_1(M)$.
In addition, $K_f$ is obviously an isotropic subspace in $H^1(M,\C)$ with respect to the cup-product. It should be maximal,
since any strictly larger isotropic subspace would contradict the fact that $K_f$ is an irreducible component of $\RR_1(M)$.

\endproof

Now we can state the following key result, which can be regarded as a strengthening of Remark \ref{Rk}.

\begin{prop} \label{prop2.5}
Let $f:M \to S$ be a surjective morphism with a generic connected fiber, such that $S=\PP^1 \setminus B$ with $\chi(S)<0$.
Then $K_f=f^*(H^1(S,\C))$ is a (rationally defined) maximal isotropic subspace in $H^1(M,\C)$ with respect to the cup-product and for any character $\rho \in \T(M)$ with $\LL_{\rho}|F=\C_F$ for a generic fiber $F$ of $f$, the translate subtorus
$$W_{f,\rho}=\rho \otimes f^*(\T(S))$$
is an irreducible component in $\V_1(M)$ such that $\dim W_{f,\rho}=-\chi(S)+1 \ge 2$.

\end{prop}

In the proof, we use the following version of {\it projection formula}, which is used very often, e.g. \cite{A},
\cite{L1}, but for which I was not able to find a reference.

\begin{lem} \label{lem2}
For any local system $\LL_1$ on $M$ and any local system $\LL_2$ on $S$, one has
$$(Rf_*\LL_1)\otimes \LL_2=Rf_*(\LL_1 \otimes f^{-1}\LL_2).$$
\end{lem}

\proof

To prove this Lemma, we start with the usual projection formula, i.e. with the above notation
\begin{equation} \label{pf1}
(Rf_!\LL_1)\otimes \LL_2=Rf_!(\LL_1 \otimes f^{-1}\LL_2)
\end{equation}
see Thm. 2.3.29, p.42 in \cite{D2}. Let $Z$ be a connected smooth complex algebraic variety of dimension $m$.
Then the dualizing sheaf $\omega _Z$ is just $\C_Z[2m]$ and $D_Z\LL=\LL^{\vee}[2m]$ for any local system $\LL$ on $Z$, see Example 3.3.8, p.69 in \cite{D2}. Note also that for two bounded constructible complexes
$\A^*$ and $\B^*$ in $D_c^b(Z,\C)$ we have the isomorphisms
\begin{equation} \label{dual1}
D_Z\A^*\otimes \B^*=RHom(\A^*,\omega _Z)\otimes \B^*=RHom(\A^*,\omega _Z\otimes \B^*)=
\end{equation}
$$=RHom(\A^*, \B^*)[2m].$$
 It follows that
\begin{equation} \label{dual2}
D_Z(\A^*\otimes \B^*)=RHom(\A^* \otimes \B^*,\omega _Z)=RHom(\A^*, RHom(\B^*, \omega _Z))=
\end{equation}
$$=D_Z\A^* \otimes D_Z \B^*[-2m].$$
For the second isomorphism here we refer to Prop. 10.23, p.175 in \cite{Bo}.
Apply now the duality functor $D_S$ to the projection formula \eqref{pf1}. In the left hand side we get
$D_S((Rf_!\LL_1)\otimes \LL_2)=D_S(Rf_!\LL_1)\otimes D_S(\LL_2)[-2]=Rf_*(D_M\LL_1)\otimes D_S(\LL_2)[-2]=
Rf_*(\LL_1^{\vee}  )\otimes \LL_2  ^{\vee}[4].$
Except the isomorphisms explained above we have used here the isomorphism $D_SRf_!=Rf_*D_M$, see Cor. 4.1.17,
p.90 in \cite{D2}. Similarly, the in the right hand side we get
$D_SRf_!(\LL_1 \otimes f^{-1}\LL_2)=Rf_*D_M(\LL_1 \otimes f^{-1}\LL_2)=Rf_*(\LL_1^{\vee} \otimes (f^{-1}\LL_2 ) ^{\vee})[4].$ Since $(f^{-1}\LL_2 ) ^{\vee}=f^{-1}(\LL_2  ^{\vee})$ and since any local system is the dual of its own dual, the proof is completed.

\endproof

\proof Now we prove Proposition \ref{prop2.5}.
The first claim follows from Remark \ref{Rk} and Corollary \ref{cor1.0}.
 Let $\LL_1=\LL_{\rho}$ and $\LL_2$ be any rank 1 local system on $S$. To estimate $\dim H^1(M, \LL_1 \otimes f^{-1}\LL_2)$, we use the Leray spectral sequence
$$E_2^{p,q}=H^p(S,R^qf_*(\LL_1 \otimes f^{-1}\LL_2))$$
converging to $H^{p+q}(M,\LL_1 \otimes f^{-1}\LL_2)$. This spectral sequence degenerates at $E_2$ since
$E_2^{p,q}=0$ for $p \notin \{0,1\}$ by Artin Theorem, see Thm.4.1.26, p.95 in \cite{D2}. By Lemma
\ref{lem2} we have 
$$R^qf_*(\LL_1 \otimes f^{-1}\LL_2)=R^qf_*(\LL_1 )\otimes \LL_2.$$
In particular, in this way, the above spectral sequence yields the following exact sequence
\begin{equation} \label{ES}
0 \to H^1(S,R^0f_*(\LL_1 )\otimes \LL_2) \to H^1(M, \LL_1 \otimes f^{-1}\LL_2) \to H^0(S,R^1f_*(\LL_1 )\otimes \LL_2) \to 0.
\end{equation}
Note that $\F=R^0f_*(\LL_1 )$ is in general no longer a local system on $S$, but a {\it constructible sheaf}.
By definition, it exists a minimal finite set $\Sigma \subset S$  such that $\F | (S \setminus \Sigma)$ is a local system. If $S' \subset S$ is a Zariski open subset such that the restriction $f':M' \to S'$
with $M'=f^{-1}(S')$, is a topologically locally trivial fibration, it follows that $\F|S'$ is a local system of rank 1.
Indeed, for $s \in S'$ we have
$$\F_s=\lim_{s\in D}\F(D)=\lim_{s\in D}H^0(f^{-1}(D),\LL_1)=\C.$$
Here the limit is taken over all the sufficiently small open discs $D$ in $S$ centered at $s$, and the last
equality comes from the fact that the inclusion $F_s=f^{-1}(s) \to f^{-1}(D)$ is a homotopy equivalence and
$\LL_1|F_s=\C_{F_s}$ (recall that $F_s$ is connected). In particular $\Sigma \subset S \setminus S'$. The above
argument shows also that $\dim \F_s \le 1$ for any $s \in \Sigma$ (since $f^{-1}(D)$ is connected as well).

To estimate $\dim H^1(S,\F \otimes \LL_2)$ we compute
$$\chi(S,\F \otimes \LL_2   )=\dim H^0(S,\F \otimes \LL_2)-\dim H^1(S,\F \otimes \LL_2)$$
using Thm. 4.1.22, p.93 in \cite{D2}. We get 
$$\chi(S,\F \otimes \LL_2)=\chi(S \setminus \Sigma)+\sum_{s\in \Sigma}\dim \F_s.$$
It follows that
\begin{equation} \label{dim1}
\dim H^1(S,\F \otimes \LL_2)=\dim H^0(S,\F \otimes \LL_2) -\chi(S)+\sum_{s\in \Sigma}(1-\dim \F_s) \ge -\chi(S).
\end{equation}
It follows that the translate torus $W_{f,\rho}=\rho \otimes f^*(\T(S))$ is contained in some irreducible component $W= \rho_1 \otimes f^*_1(\T(S_1))$ of $\V_1(M)$. Moreover, we can take $\rho=\rho_1$.
Then, taking the tangent spaces of $W_{f,\rho}$ and of $W$ at the common point $\rho$ we get
$$K_f \subset f^*_1(H^1(S_1,\C)).$$
But $f^*_1(H^1(S_1,\C))$ is an isotropic subspace, and $K_f$ is maximal with this property, hence
$K_f = f^*_1(H^1(S_1,\C))$, and therefore $W_{f,\rho}=W$.

\endproof

\begin{rk} \label{rk0.15}
In the proof of Proposition 1.7 in Arapura \cite{A}, section V, it is shown that
for any $\LL_1 \in \T(M)$ with trivial  restriction to generic fibers of $f$, one has 
$$ H^0(S,R^1f_*(\LL_1 )\otimes \LL_2)= 0$$
for all but finitely many local systems $\LL_2 \in \T(S)$. Actually the proof in Arapura \cite{A}  is given in the case
$\LL_1=\C_M$, but the same approach using properties of constructible sheaves yields the general case stated above.

 Using this, the exact sequence \eqref{ES} yields
$$ H^1(M, \LL_1 \otimes f^{-1}\LL_2)=H^1(S,R^0f_*(\LL_1 )\otimes \LL_2)$$
for all but finitely many local systems $\LL_2 \in \T(S)$. In particular, when $\LL_1=\C_M$, then
$R^0f_*(\LL_1 )=\C_S$ and the formula \eqref{dim1} implies that
$$\dim H^1(M,  f^{-1}\LL_2)=-\chi(S)$$
for all but finitely many local systems $\LL_2 \in \T(S)$.

\end{rk}

\begin{cor} \label{cor1.1}
If $f:M \to \C^*$ is a surjective morphism with connected generic fiber $F$ and $\rho \in \T(M)$ is a torsion character such that  $W=\rho \otimes f^*(\T(\C^*))$ is an irreducible component of $\V_1(M)$ with $\dim W =1$, then
$K_f=f^*(H^1(\C^*,\C))$ is a rationally defined  maximal isotropic subspace in $H^1(M,\C)$ with respect to the cup-product.

Conversely, let $f:M \to \C^*$ be a morphism with a generic connected fiber.
Assume that $K_f=f^*(H^1(\C^*,\C))$ is a rationally defined maximal isotropic subspace in $H^1(M,\C)$ with respect to the cup-product.
Then for any character $\rho \in \T(M)$ with $\LL_{\rho}|F=\C_F$ for a generic fiber $F$ of $f$ and $\rho \notin f^*(\T(\C^*))$ , the translated subtorus
$$W_{f,\rho}=\rho \otimes f^*(\T(\C^*))$$
is either an irreducible component in $\V_1(M)$ such that $\dim W_{f,\rho}=1$, or $ H^1(M,\LL)=0$ for 
$\LL \in W_{f,\rho}$ with finitely many exceptions.
\end{cor}

\proof Assume that $K_f$ is not maximal, and let $K \supset K_f$ be a maximal isotropic subspace in $H^1(M,\C)$ with respect to the cup-product. Then $\dim K \ge 2$ and there is a morphism $f_1:M \to S_1$  surjective, with connected generic fiber such that $K=K_{f_1}$. But then $W$ is strictly contained in
 $W_1=\rho \otimes f_1^*(\T(S_1))$, which is a component of  $\V_1(M)$ by Proposition  \ref{prop2.5},   a contradiction.
 
 For the converse part, just note that, using the same argument as above, $W_{f,\rho}$ cannot be strictly contained in a component of $\V_1(M)$.

\endproof

\begin{rk} \label{rk0.2}

(i) Unlike the case of maximal isotropic subspaces in $H^1(M,\C)$ of dimension at least two which are automatically rationally defined, see Remark \ref{rk0.1}, there are a lot of non  rationally defined
 maximal isotropic subspaces in $H^1(M,\C)$ of dimension 1 as soon as a rationally defined one exists. To see this, use the semi-continuity of the dimension of $H^1(H^*(M,\C),\mu_z)$ with respect to $z \in H^1(M,\C)$.

(ii) If $E$ is a maximal isotropic subspace in $H^1(M,\C)$ of dimension at least two, then it has at least one associated component $W_E$ in  $\V_1(M)$ corresponding to $E$  under the bijection in Theorem \ref{thm1}.
On the other hand, if  $E$ is a rationally defined maximal isotropic subspace in $H^1(M,\C)$ of dimension
1, it is quite possible that there is no associated component in  $\V_1(M)$.
As an explicit example, consider the case of the line arrangement $xyz=0$ in $\PP^2$.
Then $M=(\C^*)^2$ and any 1-dimensional subspace  $E_{\Q}\subset H^1(M,\Q)$ gives rise to a
 rationally defined maximal isotropic subspace in $H^1(M,\C)$ of dimension 1. However, it is well known that
 $\V_1(M)=\{1\}$ in this case.

\end{rk}

We say that an irreducible component $E$ of some $\RR_m(M)$ (which is a vector subspace in $ H^1(M)$) is a {\it coordinate component},
resp. a {\it non-essential component}, if it corresponds under the above isomorphism to a
coordinate (resp. non-essential) component of $\V_m(M)$. Proposition \ref{prop3} can be reformulated as follows.

\begin{prop} \label{prop4}
An irreducible component of $\RR_m(M)$ is non-essential if and only if it is a coordinate component.
\end{prop}

An  irreducible component $E$ of some $\RR_m(M)$ which is not a coordinate component is called
a  {\it global component} in  \cite{FY}. This is compatible with our Definition  \ref{def1} above.

\medskip

Let us consider the  component $W_E$ introduced in Remark  \ref{rk0.2}, in the case $\dim E \ge 2$.
Then $W_E$ corresponds to a mapping $f_E:M \to S_E$ as in Theorem \ref{thm1.5}. We have the following result.

\begin{prop} \label{prop4.5}
If the arrangement $C$ is given, then the mapping $f_E$ is determined by the vector subspace $E \subset H^1(M,\C)$. 

\end{prop}

\proof

First we know that $S_E$ is obtained from $\PP^1$ by deleting a subset $B$, with $|B|-1=\dim E$.
Let $k:=\dim E$ and assume that the points in $B$ are $(0:1)$ and $(1:b_j)$, for some $b_j \in \C$,
$j=1,...,k$. If $(u:v)$ are the homogeneous coordinates on $\PP^1$, then the cohomology group $H^1(S_E,\C)$ has a basis given by
$$\omega _j=\frac{d(v-b_ju)}{v-b_ju}-\frac{du}{u}$$
where $j=1,...,k.$ As explained in Proposition \ref{prop2}, the mapping $f_E$ corresponds to a pencil
$(P,Q)$, where $P$ and $Q$ are homogeneous polynomials in $\C[X,Y,Z]$, of the same degree and without common factors. In terms of this pencil, one has
$$f_E^*(\omega _j)=\frac{d(Q-b_jP)}{Q-b_jP}-\frac{dP}{P}$$
So, in down-to-earth terms, the question is: how to determine the pencil $(P,Q)$ from the vector space
of 1-forms with logarithmic poles
$$E=<f_E^*(\omega _j)~~|~~j=1,...,k>?$$
Using only logarithmic poles allows us to work with rational differential forms rather than cohomology classes
and this is essential for this proof.
Start with the curve $C_1$ in the curve arrangement $C$ and consider the subset
$$E_1=\{\omega \in E~~|~~ \int _{\gamma_1}\omega = 0\}.$$
Two cases may occur.

\medskip

\noindent Case 1. ($E_1= E$) This case occur exactly when $C_1$ is not a connected component in any of the
$(k+1)$ special fibers of the pencil $(P,Q)$ corresponding to the set $B$. If this happens, we discard the curve $C_1$ and test the next curve $C_2$ and so on.

\medskip

\noindent Case 2. ($E_1 \ne E$) Then $C_1$ is a component of a special fiber of the pencil, say of the fiber
$\CC_{b_1}: Q-b_1P=0$. Note that any form $\omega \in E$ can be written as a sum
$$\omega = \sum_ja_jf_E^*(\omega _j)$$
and, with this notation, one has
$$\int _{\gamma_1}\omega=2\pi i a_1 m(C_1)$$
where $i^2=-1$ and $m(C_1)$ is the multiplicity of the irreducible curve $C_1$ in the divisor $\CC_{b_1}$.
It follows that $\omega \notin E_1$ if and only if $a_1 \ne 0$. This shows that
$Q-b_1P$ is the G.C.D. of the denominators of the forms $\omega \in E \setminus E_1$.
After we have determined the polynomial $Q-b_1P$ as above, we discard all the curves $C_j$ in the arrangement
which are contained in the support of the divisor $\CC_{b_1}$. From the remaining curves in $C$, we can find a new curve, say $C_2$, such that
$$E_2:=\{\omega \in E~~|~~ \int _{\gamma_2}\omega = 0\} \ne E.$$
(such a curve exists since $|B|\ge 3$). Then $C_2$ is a component in a new fiber of the pencil $(P,Q)$,
say of the fiber
$\CC_{b_2}: Q-b_2P=0$ and hence $Q-b_2P$ is the G.C.D. of the denominators of the forms $\omega \in E \setminus E_2$. The two homogeneous polynomials $Q-b_1P$ and $Q-b_2P$ span the same vector space as the polynomials $P$ and $Q$, i.e. they determine the same pencil up to an automorphism of $\PP^1$.

\endproof

Now we treat the special case of rationally defined maximal isotropic subspace in $H^1(M,\C)$ of dimension 1.

\begin{prop} \label{prop4.6}

Let $E$ be a rationally defined maximal isotropic subspace in $H^1(M,\C)$ of dimension 1. Then there is a
surjective mapping $f_E:M \to \C^*$ with connected generic fiber such that $E=f_E^*(H^1(\C^*,\C))$.
\end{prop}

\proof

Let $f_j=0$ be a homogeneous reduced equation for the component $C_j$ in the curve arrangement $C$, for $j=1,...,r.$ Assume that $\deg f_j=d_j$. A basis of the $(r-1)$-dimensional vector space  $H^1(M,\Q)$
is given by the 1-forms
$$\eta_k=\frac{1}{d_k} \frac{df_k}{f_k}-\frac{1}{d_{k+1}} \frac{df_{k+1}}{f_{k+1}}$$
for $k=1,...,r-1.$ Using this, we see that any 1-dimensional subspace in  $H^1(M,\Q)$ has a unique
generator (up to a $\pm $-sign) of the form
$$\eta=\sum _{j=1,r}m_j\frac{df_j}{f_j}$$
where $m_j$ are relatively prime integers, i.e. $G.C.D.(m_1,...,m_r)=1$, such that
$$\sum _{j=1,r}d_jm_j=0.$$
It follows that the rational fraction
$$f=\prod _{j=1,r}f_j^{m_j}$$
is homogeneous of degree $0$ and hence induces a morphism $f:M \to \C^*$. The fact that $f$ is surjective is
obvious, while the connectivity of the generic fiber follows from Bertini's Theorem, see \cite{Sch}, p. 79, using the condition $G.C.D.(m_1,...,m_r)=1$. The equality $E=f^*(H^1(\C^*,\C))$ is obvious by taking
$\frac{dt}{t}$ as a basis of $H^1(\C^*,\C)$.
\endproof

\section{Minimal arrangements} \label{sec:four}

In this section we prove the following result which applies to an arbitrary
plane curve arrangement.

\begin{thm} \label{thm2}
Let $C=\cup _{i=1,r}C_i$ be a plane curve arrangement in $\PP^2$, having $r$ irreducible components $C_i$, for $i=1,r$. Let $M=\PP^2 \setminus C$ be the corresponding complement.
Then, for $d \ge 2$, the following are equivalent.

\medskip

\noindent (i) there is a global
$d$-dimensional irreducible component $E$ in the resonance variety  $\RR_1(M)$;

\medskip

\noindent (ii) there is a global 
$d$-dimensional irreducible component $W$ in the characteristic variety  $\V_1(M)$;

\medskip

\noindent (iii) there is
a pencil $\CC$ of plane curves on $\PP^2$ with an irreducible generic member and having
 $d+1$ fibers $\CC_t$ whose reduced supports form a partition
of the set of irreducible components  $C_i$, for $i=1,r$.

\medskip

\noindent Moreover, for $d=1$, (i) always fails and (ii) implies (iii).

\end{thm}

\proof

For $d \ge 2$, the equivalence between $(i)$ and $(ii)$ follows from Theorem \ref{thm1} and  
Theorem \ref{thm1.5}, claim (vi).
And the equivalence between $(ii)$ and $(iii)$ follows from Arapura's results recalled in subsection \ref{sec:3.3} combined with Proposition \ref{prop2} and Corollary  \ref{cor1}. The case $d=1$ follows from
Theorem \ref{thm1.5}, claim (i).

\endproof

Note that the condition $(ii)$ above can be reformulated as the existence of a global  irreducible component $W$ in the characteristic variety  $\V_1(M)$ such that $1\in W$ and $\dim H^1(M,\LL)=d-1 \ge 1$ for a generic local system
$\LL \in W$.

\medskip
 
The following result is similar to Theorem 4.1.1 in  Libgober \cite{L1} and closely related to the discussion in  \cite{LY}, just before Proposition 7.2.

\begin{cor} \label{cor2}
Assume the equivalent statements in Theorem \ref{thm2} above hold. Let $f:\PP^2 \to \PP^1$ be the rational morphism associated to the pencil  $\CC$. Let
$\pi: X \to \PP^2$ be a sequence of blowing-ups such that $g=f \circ \pi$
is a regular morphism on $X$. 
If $\tilde C$ denotes the proper transform of the (reduced) curve $C$ under $\pi$, then the self intersection number of  $\tilde C$ is non-positive, i.e.
$$\tilde C \cdot \tilde C \le 0.$$

\end{cor}

\proof

There is a partition of $\tilde C =\cup _{i=1,d+1} ~ \tilde C_i$ of $\tilde C$ as a union of disjoint curves $\tilde C_i$, such that $g(\tilde C_i)=b_i$ for $i=1,...,d+1$.
It follows that 
$$ \tilde C \cdot \tilde C =\sum  _{i=1,d+1}\tilde C_i \cdot \tilde C_i.$$
Each curve is contained in the support of the positive divisor $g^{-1}(b_i)$,
and hence by Zariski's Lemma, see  \cite{BPV}, p. 90, we get
$$\tilde C_i \cdot \tilde C_i \le 0.$$

\endproof

\begin{ex} \label{ex1}
Assume that $C$ is a line arrangement, i.e. $d_j=1$ for all $j=1,...,r$.
Let $\X$ be the base locus of a pencil as in Corollary \ref{cor2}, $k=d+1$ the number of sets in the associated partition $( \tilde C_i)_i$  of the set of lines in $C$.
Then it is shown in \cite{FY} that the following hold.

\medskip

\noindent (i) For each base point $p \in \X$, the multiplicity
$$ n_p=\mult_p( \CC_b)$$
is independent of $b \in B$.

\medskip

\noindent (ii) $\sum _{p\in \X}n_p=D^2$, where $D$ is the degree of the pencil.

\medskip

\noindent (iii)  $\sum _{j=1,r}m(C_j)=kD$, where $m(C_j)\ge 1$ is the multiplicity with which $C_j$ occurs in the corresponding fiber of the pencil.

\medskip

To resolve the indeterminacy points of the associated pencil (i.e. to determine the map $\pi:X \to \PP^2$),  one has in this case just to blow-up once the points in the base locus $\X$. This is a direct consequence of the property (i)
above.

Assume moreover that  $m(C_j)= 1$ for all $j=1,r$, i.e. all $\CC_b$ for $b \in B$ are reduced.

Then, again by (i) above, it follows that $\mult _p(C)=kn_p$ for any $p \in \X$.
On the other hand, by (iii) we get $\deg C = kD$.
Finally, in this very special  case we get 
$$\tilde C \cdot \tilde C= C \cdot C -\sum_{p\in \X} \mult _p(C)^2 = k^2D^2-
 k^2D^2=0$$
since $ C \cdot C = (\deg C)^2$. This happens for instance in Example 3.4 in
\cite{FY}: the Ceva arrangement given by the pencil $ax^d+by^d +cz^d=0$ with
$(a:b:c) \in \PP^2$ satisfying $a+b+c=0$. There are 3 special fibers, corresponding to $x^d-y^d$, $y^d-z^d$, $z^d-x^d$.

 In the case of a general line arrangement, the condition $\tilde C \cdot \tilde C\le 0$ may bring new non-trivial information on the arrangement. In particular it can be used as a test for candidates to the base locus $\X$ of a pencil associated to a given arrangement. In  Example 3.6 in
\cite{FY} the $B_3$-arrangement consists of 9 lines, and the base locus  $\X$
 consists of 3 points of multiplicity 4 and 4 other points of multiplicity 3.
As a result we have
$$\tilde C \cdot \tilde C= C \cdot C -\sum_{p\in \X} \mult _p(C)^2 = 81-3\times 16 -4\times 9=-3.$$
This latter arrangement is associated to the pencil $((x^2-y^2)z^2:(y^2-z^2)x^2)$ which has again  3 special fibers (this time non-reduced!), corresponding to $(x^2-y^2)z^2   $, $(y^2-z^2)x^2  $, $(x^2-z^2)y^2$.
Up to a linear change of coordinates, this is the same $B_3$-arrangement as in the final subsection (6.7) below.

\end{ex}

\section{Fibered complements and $K(\pi,1)$-spaces} \label{sec5}

Let $f:M \to S$ be a morphism associated to a plane curve pencil $\CC$
with base locus $\X$, as in
Theorem  \ref{thm2} and Corollary  \ref{cor2} above.
Consider a fiber $\CC_s$ of the  pencil $\CC$ corresponding to 
$s\in S=\PP^1 \setminus B$. We say that  $\CC_s$ is a {\it special fiber} of  $\CC$ if either
$\CC_s \setminus \X $ is singular , or if $\CC_s \setminus \X $ is smooth
and exists a point $p \in \X$ such that
$$ \mu(\CC_s,p) > \min_{t \in S}\mu(\CC_t,p),$$
where $\mu$ denotes the Milnor number of an isolated singularity.

Let $\CC_{spec}$ be the union of all the special fibers in the pencil  $\CC$.
(There is a finite number of such fibers, and they are easy to identify, see
\cite{LW},  \cite{DM}). Let $B'$ be the union of $B$ and the set of all $s \in S$ such that $\CC_s$ is a special fiber. We call $C' = C\cup \CC_{spec}$ the {\it extended plane curve arrangement} associated to the  plane curve arrangement
$C$ and denote by $M'$ the corresponding complement. We set $S'=\PP^1 \setminus B'$.
With this notation, we  have the following result.

\begin{prop} \label{prop5}
The restriction $f':M' \to S'$ is a locally topologically trivial fibration.
The plane curve arrangement  complement $M'$ is a $K(\pi,1)$-space.

\end{prop}

\proof

Note that for any $p \in \X$, the family of plane curve isolated singularitis
$(\CC_s,p)$ for $s \in S'$ is a $\mu$-constant family. Using the relation between  $\mu$-constant families and Whitney regular stratifications, as well as 
Thom's First Isotopy Lemma, see for instance  \cite{D1}, pp. 11-16 and especially the proof of Proposition (1.4.1) on p. 20,
we get the first claim above. The second claim is an obvious consequence, as explained already in \cite{FY}.

\endproof

\begin{rk} \label{rk3}
With the notation from Theorem \ref{thm1.5}, Libgober has remarked in the proof of Lemma 1.4.3 in
 \cite{L1} that
$$\dim H^1(M,\LL_{\rho _W}\otimes f^*_W(\LL))=\dim H^1(S_W,f_{W*}(\LL_{\rho _W})\otimes \LL)$$
for almost all local systems $\LL \in \T(S_W)$ (this comes from the exact sequence  \eqref{ES} by proving the vanishing of the last term). This implies that a morphism $f_W:M \to S_W$ as in 
 Theorem \ref{thm1.5} is never a locally trivial fibration in the case $S_W=\C^*$ (this is related to the case $\Sigma =\emptyset$ in formula  \eqref{dim1}).
By definition, $\dim H^1(M,\LL_{\rho _W}\otimes f^*_W(\LL)) \ne 0$ while the second dimension would vanish in such a situation. Indeed then $f_{W*}(\LL_{\rho _W})$ would be a local system on $\C^*$,
given by an automorphism $A$ of some $\C^n$ and we can take $\LL$ to have a monodromy different from the eigenvalues of $A$.

\end{rk}

\begin{ex} \label{ex2}

We discuss now the following curve arrangement, considered already in 
Example 4.8 in \cite{DiM}. The curve $C$ consists of the following:
three lines $C_1: ~x=0$, $C_2: ~y=0$ and $C_3: ~z=0$ and a conic $C_4: ~x^2-yz=0$. The corresponding pencil can be chosen to be $f=(x^2:yz)$, and  the set
$B=\{b_1,b_2,b_3\}$ is given by $b_1=(0:1)$, $\CC_1=2C_1$,
 $b_2=(1:0)$, $\CC_2=C_2 \cup C_3$, and  $b_3=(1:-1)$, $\CC_3=C_4$.

The base locus of this pencil is $\X=\{p_1=(0:0:1), p_2=(0:1:0)\}$ and it is easy to check that there no special fibers. It follows that $f:M \to S$ is a locally topologically trivial fibration with fiber
$\C^*$ (a smooth conic minus two points). Hence  the complement $M$ of this 
curve arrangement is a $K(\pi,1)$-space, where the group $\pi$ fits into an exact sequence
$$ 1 \to \Z \to \pi \to \FF_2\to 1$$
with $ \FF_2$ denoting the free group on two generators.

\medskip

Note also that 
$$\mult_{p_1}\CC_1=2 > 1=\mult_{p_1}\CC_2,$$
hence the property (i) in Example \ref{ex1} does not hold for arbitrary curve arrangements.

\medskip

A local system $\LL$  on $S$ is given by a triple $(\lambda_1,\lambda_2, \lambda_3)$
with $\lambda_1 \cdot \lambda_2 \cdot \lambda_3=1$. With this notation, the pull-back local system $f^*\LL$ is given by $(\rho_1, \rho_2,\rho_3, \rho_4)$,
where $\rho_1= \lambda_1^2$,  $\rho_2=\rho_3= \lambda_2$ and 
 $\rho_4= \lambda_3$ (recall the discussion before Corollary \ref{cor1}). With this, it is easy to check that
 the irreducible component given by $f^*\T(S)$ coincide with the irreducible component obtained in 
Example 4.8 in \cite{DiM} by computation using the associated integrable connections.

Indeed, the $\rho_i$'s satisfy the equation $\rho_1 \cdot \rho_2  \cdot \rho_3 \cdot \rho_4^2=1$,
see Equation \eqref{eq5}, and hence we can use  $(\rho_1, \rho_2, \rho_4)$ to parametrize the torus $\T(M)$.
It follows that  the irreducible component given by $W=f^*\T(S)$ is parametrized by
 $\rho_1= \lambda_1^2$,   $\rho_2= \lambda_2$ and $\rho_4= \lambda_3$, i.e. $W$ is given by the equation
$$\rho_1\cdot \rho_2 ^2  \cdot \rho_4^2=1$$
which appears in Example 4.8 in \cite{DiM} with slightly different notation.

\medskip

Let $\wC$ be the proper transform on $C$ under the blowing-up $\widehat \pi: \widehat X \to \PP^2$
of the two points in $\X$. Then
$$\wC \cdot \wC=C\cdot C- 3^2- 3^2=7>0.$$
This is not in contradiction with Corollary  \ref{cor2}, since in order to resolve the indeterminacy points of $f$ in this case we have to blow the points $\widehat p_1, \widehat p_2 \in \widehat X$ corresponding to the tangents of the conic $C_4$ at $p_1$ and $p_2$. The new multiplicities are 
$$\mult_{\widehat p_1}\wC = \mult_{\widehat p_2}\wC=3$$
and hence for the new proper transform $\tilde C$ we get

$$\tilde C \cdot \tilde C= \wC \cdot \wC - 3^2- 3^2=-11.$$

\end{ex}

\section{On the translated components of $\V_1(M)$} 

Let $W$ be a translated irreducible component  of $\V_1(M)$, i.e. $1 \notin W$. Then, as in Theorem  \ref{thm1.5}, there is a torsion character  $\rho \in \T(M)$ and a surjective morphism $f:M\to S$ with connected generic fiber $F$ such that
\begin{equation} \label{k1}
 W=\rho f^*(\T(S))
\end{equation}
We say in this situation that the component $W$ is associated to the mapping $f$.
In this section we give detailed information on the torsion character $\rho \in \T(M)$ in terms of the geometry of the associated
mapping $f:M \to S$.

\subsection{The general setting} 

Let $F$ be the generic fiber of the mapping $f:M \to S$, i.e. $F$ is the fiber of the topologically locally trivial fibration
$f':M' \to S'$ associated to $f$ as in the previous section. Then, we have an exact sequence
\begin{equation} \label{es1}
H_1(F)   \stackrel{i'_*} \longrightarrow     H_1(M')  \stackrel{f'_*} \longrightarrow      H_1(S') \to 0 
\end{equation}
as well as a sequence 
\begin{equation} \label{es2}
H_1(F)  \stackrel{i_*} \longrightarrow      H_1(M)   \stackrel{f_*} \longrightarrow      H_1(S) \to 0 
\end{equation}
which is not necessarily exact in the middle, i.e. the group
\begin{equation} \label{t}
T(f)=\frac{\ker f_*}{\im i_*}
\end{equation}
is in general non-trivial. Here $i:F\to M$ and $i':F \to M'$ denote the inclusions, and homology is taken with $\Z$-coefficients
if not stated otherwise. The group $T(f)$ is a finite abelian group according to Theorem  \ref{thm5}, an apparently not obvious fact.

This group was studied in a compact (proper) setting by Serrano, see \cite{Se},
but no relation to local systems was considered there. On the other hand, this compact situation was also studied by A. Beauville in \cite{Beau}, with essentially the same aims as ours. However, the actual results 
 are distinct, because of the key role played in our case by the partially deleted fibers, corresponding to the two types of terms in the sum in formula \eqref{f2} below.

The second sequence induces an obvious  exact sequence
\begin{equation} \label{es3}
0 \to T(f) \to \frac{H_1(M)}{\im i_*} \stackrel{f_*} \longrightarrow      H_1(S) \to 0. 
\end{equation}
Since $ H_1(S)$ is a free $\Z$-module, applying the fuctor $\Hom(-,\C^*)$ to the  exact sequence \eqref{es3}, we get a new exact sequence
\begin{equation} \label{es4}
1 \to \T(S) \to \T(M)_F \to \Hom(T(f),\C^*) \to 1.
\end{equation}
Here $\T(M)_F$ is the subgroup in $\T(M)$ formed by all character $\chi: H_1(M) \to \C^*$ such that
$\chi \circ i_*=0$. This means exactly that the associated local system $\LL _{\chi}$ by restriction to $F$
yields the trivial local system $\C_F$.

The torsion character $\rho \in \T(M)$ which occurs in  \ref{k1} is in this subgroup  $\T(M)_F$, see Theorem  \ref{thm1.5}, (ii).
Moreover, this  character $\rho$ is not unique, but its class 
\begin{equation} \label{k2}
{\tilde \rho} \in \frac{\T(M)_F}{\T(S)} \simeq  \Hom(T(f),\C^*)
\end{equation}
is uniquely determined. From now on, we will regard ${\tilde \rho} \in  \Hom(T(f),\C^*)$. Hence, to understand the possible choices for
 ${\tilde \rho}$, we have to study the group $T(f)$.

\subsection{The computation of the group $T(f)$} 

In order to simplify the presentation, we assume in this subsection that at least one of the curves $C_j$ in  the curve arrangement
$C$ is a line. This covers the case of line arrangements and of curve arrangements in the affine plane $\C^2$.
More specifically, we assume that $C_1$ is the line at infinity in $\PP^2$ and hence $M =\C^2  \setminus \cup_j(\C^2 \cap C_j)$.
We assume that $\infty \in B$ and set $B_1= B \setminus \{\infty\} \subset \C$.

 For each $b \in B$, consider the following divisor
on $\C^2$
\begin{equation} \label{f1}
D_b= \overline {g_1^{-1}(b)} \cap \C^2=\sum _a m_{ba}C_{ba}
\end{equation}
where $g_1$ is the extension of $f$ from the proof of Proposition  \ref{prop2}, $ m_{ba} \ge 1$ are integers and $C_{ba}$
are irreducible curves in the arrangement $C$.

Recall the larger set $B' \supset B$ obtained from $B$ by adding the bifurcation points of $f:M \to S$. Set $C(f)=B' \setminus B$ and assume from now on that $C(f)$ is nonempty. (Otherwise $f$ is a fibration and hence $T(f)=1$).
 For each $c \in C(f)$, consider the following divisor
on $\C^2$
\begin{equation} \label{f2}
D_c= \overline {g_1^{-1}(c)} \cap \C^2  =\sum _{a'} m_{ca'}'C_{ca'}'+\sum _{a''} m_{ca''}''C_{ca''}''
\end{equation}
where $g_1$ is the extension of $f$ as above, $ m_{ca'}'\ge 1$ and $ m_{ca''}''\ge 1$ are integers,
the irreducible curves $C_{ca'}'$ are curves in our arrangement (corresponding to case (2) in Proposition  \ref{prop2},
and $C_{ca''}''$ are the new curves to be deleted from $M$ in order to obtain $M'$. Since $f:M \to S$ is a surjection,
it follows that there is at least one term in the second sum in the equality  \eqref{f2}. At least some of the first type sums are non-trivial if and only if the arrangement $C$ is special.

Note that $M$ is obtained from  $\C^2$ by deleting the curves $C_{ba}$, $C_{ca'}'$ and possibly some horizontal
components $C_h$. For each irreducible curve $Z$ in  $\C^2$ we denote by $\gamma(Z)$ the elementary oriented loop
associated to $Z$. It follows that $H_1(M)$ is a free $\Z$-module with a basis given by
\begin{equation} \label{base1}
\gamma(C_{ba}), \gamma(C_{ca'}'),      \gamma(C_h)
\end{equation}

Similarly,  $H_1(S)$ is a free $\Z$-module with a basis given by $\delta_b$ for $b \in B_1$, where $\delta_b$
is an elementary loop based at $b$ as subsection (3.1). In term of these bases, the morphism $f_*:H_1(M) \to H_1(S)$
is described as follows.
Let $\alpha \in H_1(M)$ be given by
\begin{equation} \label{m1}
\alpha = \sum \alpha _{ba}   \gamma(C_{ba})   +  \sum  \alpha _{ca'}   \gamma(C_{ca'}')   + \sum \alpha_h \gamma(C_h).
\end{equation}
Here $a \in A_b$, an index set depending on $b$. This is not written explicitely, in order to keep the notation simpler. Similar remarks apply to the indices $a'$ and $a''$ below.
Then
\begin{equation} \label{m2}
f_*(\alpha) = \sum _{b \in B_1}( \sum _a  m_{ba} \alpha _{ba}  -  \sum _a m_{\infty a}   \alpha _{\infty a})\delta_b  .
\end{equation}
In particular $\alpha \in \ker f_*$ if and only if
\begin{equation} \label{m3}
 \sum _a m_{ba} \alpha _{ba}  =  \sum _a m_{\infty a} \alpha _{\infty a}
\end{equation}
for all $b \in B_1$. Next, $M'$ is obtained from $M$ by deleting the curves $C_{ca''}''$. Hence, it follows that $H_1(M')$ is a free $\Z$-module with a basis given by
\begin{equation} \label{base2}
\gamma(C_{ba}), \gamma(C_{ca'}'),      \gamma(C_h),   \gamma(C_{ca''}'').
\end{equation}
The inclusion $j: M' \to M$ induces a morphism $j_*:H_1(M') \to H_1(M)$ which, at coordinate level, is just the obvious projection.
Similarly,  $H_1(S')$ is a free $\Z$-module with a basis given by $\delta_b$ for $b \in B_1$, and $\delta_c$ for  $c \in C(f)$.

 In term of these bases, the morphism $f_*':H_1(M') \to H_1(S')$
is described as follows.
Let $\beta \in H_1(M')$ be given by
\begin{equation} \label{m4}
\beta = \sum \beta _{ba}   \gamma(C_{ba})   +  \sum  \beta _{ca'}   \gamma(C_{ca'}')   + \sum \beta_h \gamma(C_h) +  \sum  \beta _{ca''}   \gamma(C_{ca''}'')  .
\end{equation}
Then $f'_*(\beta)=E_1+E_2$, where 
\begin{equation} \label{m5}
E_1 = \sum _{b \in B_1}( \sum _a  m_{ba} \beta _{ba}  -  \sum _a m_{\infty a}   \beta _{\infty a})\delta_b   .
\end{equation}
and
\begin{equation} \label{m6}
E_2= \sum_{c\in C(f)}(
\sum_{a'} m'_ {ca'}  \beta _{ca'}+ \sum_{a''} m''_ {ca''}  \beta _{ca''} -  \sum _a m_{\infty a}   \beta _{\infty a})\delta_c    .
\end{equation}
The exact sequence \eqref{es1} yields $\ker i'_*=\ker f'_*$. On the other hand $\im i_*=j_*(\im i'_*)$. It follows that 
 $\alpha \in H_1(M)$ as above is in $\im i_*$ if and only if there is a  $\beta \in H_1(M')$ such that

\medskip

\noindent (i) $\beta$ is a lifting of $\alpha$, i.e. $j_*(\beta)=\alpha$;

\medskip

\noindent (ii) $f'_*(\beta)=0.$

\medskip

In terms of our bases this means that 

\medskip

\noindent (i') $\beta _{ba} = \alpha _{ba}$,  $\beta _{ca'} = \alpha _{ca'}$,  $\beta_h = \alpha_h$; 

\medskip

\noindent (ii') The coordinates of  $\alpha$ satisfy the equation  \eqref{m3} and there is a choice of coordinates $ \beta _{ca''}$
such that one has
\begin{equation} \label{m7}
 \sum_{a''} m''_ {ca''}  \beta _{ca''} =  \sum _a m_{\infty a}   \alpha _{\infty a}-\sum_{a'} m'_ {ca'}  \alpha _{ca'}
\end{equation}
for each $c \in C(f)$.
Let $m''(c)=G.C.D.\{ m''_ {ca''} \}$ where $a''$ takes all the possible values (this set of indices is nonempty). The above equation has a solution
if and only if the right hand side is divisible by $m''(c)$.
This explains the following construction. Let $G(f) =\oplus_{c \in C(f)} \Z/m''(c)\Z$ and let
$$\theta: \ker f_* \to G(f)$$
be the  morphism sending  $\alpha \in \ker f_*$ to the element in $G(f)$ having as its $c$-coordinate the class of
$$ \sum _a m_{\infty a}   \alpha _{\infty a}-\sum_{a'} m'_ {ca'}  \alpha _{ca'}$$
modulo $m''(c)$. Let  $m'(c)=G.C.D.\{ m'_ {ca'} \}$ where $a'$ takes all the possible values (this set can be  empty  and in this case we set $m'(c)=0$).
The above discussion is summarized in the following.
\begin{thm} \label{thm5}
With the above notation, there is an isomorphism $T(f) \simeq \im \theta$. In particular, if for all $c \in C(f)$,
one has $GCD(m'(c),m''(c))=1$, i.e. none of the fibers $D_c$ is multiple, then 
$$T(f) \simeq G(f).$$
Moreover, if all the fibers  $D_c$ are reduced (i.e. all components occur with multiplicity 1), then $T(f)=1$,
and hence there are no translated components in $\V_1(M)$ associated to $f$  in this case.
\end{thm}

\begin{rk} \label{rk4}
The horizontal components $C_h$ play no role in the computation of $T(f)$, since they are obviously
in $\im i_*$.

\end{rk}

The following consequence does not rule out the possibility of a {\it translated global component}, but explains in a sense why they should be quite exceptional. Recall that by Theorem  \ref{thm2}, a translated global component corresponds to a minimal arrangement. Conversely, a translated coordinate component is more likely to occur, and then it is related to a special arrangement, as in the deleted $B_3$-arrangement revisited below.

\begin{cor} \label{cor4}
Assume that $C$ is a minimal arrangement with respect to $f:M \to S$. For $b \in B$, set $m(b)=G.C.D.\{ m_ {ba} \}$
and then $m(f)=L.C.M. \{ m(b)~~|~~ b \in B \}$. Then $T(f)$ is isomorphic to the cyclic subgroup in $G(f)$ spanned by 
$$(m(f),m(f),...,m(f)).$$
In particular, if $m''(c)|m(f)$ for all $c \in C(f)$, then $T(f)=1$.

\end{cor}

\proof
If  $C$ is a minimal arrangement, then there are no indices of type $a'$, so the morphism $\theta$ has a simpler form.
The fact that $\im \theta$ is the cyclic subgroup in $G(f)$ spanned by 
$(m(f),m(f),...,m(f))$ follows from the equations  \eqref{m3}.

\endproof

\begin{cor} \label{cor5}
If $f:M \to S$ has no multiple fibers, then there are no translated components in $\V_1(M)$
associated to $f$.
\end{cor}
The following result clarifies to a certain extent the case of 1-dimensional translate components.

\begin{prop} \label{prop6}
Let $E \subset H^1(M,\C)$ be a maximal isotropic subspace with respect to the cup product such that
$\dim E=1$. Let $f:M \to \C^*$ be the corresponding surjective morphism, with connected generic fiber $F$,
i.e. $E=f^*(H^1(\C^*))$. Assume that $m'(c)=1$ for all $c \in C$.
Then, for any nontrivial character $\tilde \rho: T(f) \to \C^*$, the associated 1-dimensional translated subtorus
$$W_{f,\rho}=\rho \otimes f^*(\T(\C^*))$$
is a component in  $\V_1(M)$.

\end{prop}

\proof
Since  $m'(c)=1$ for all $c \in C$, we obtain a system of generators for the group $T(f)$ by taking
$\theta (\gamma(C'_{ca'}))$, for  all $c \in C$ and all possible values of $a'$. Since  $\tilde \rho$
is non-trivial, it exists at least one point  $c \in C$ and one value for $a'$ such that
$$\rho(\gamma(C'_{ca'}))=\tilde \rho (\theta (\gamma(C'_{ca'}))) \ne 1.$$
If $D$ is any small disc containing $c$, it follows that
$$H^0(f^{-1}(D),\LL_{\rho})=0.$$
In the notation of the proof of Proposition \ref{prop2.5}, we get $\F_c=0$, and hence $c \in \Sigma$.
The formula  \eqref{dim1} then implies that 
$\dim H^1(S,\F \otimes \LL_2) \ge 1$ for all $ \LL_2 \in \T(\C^*)$. We conclude by applying Corollary
 \ref{cor1.1}.

\endproof
The same proof as above yields the following result, to be compared with Theorem \ref{thm1.5}, (iv).

\begin{cor} \label{cor6}
Let $f:M \to S$ be a surjective morphism, with connected generic fiber $F$, such that $\chi(S)<0$.
 Assume that $m'(c)=1$ for all $c \in C$.
Then, for any nontrivial character $\tilde \rho: T(f) \to \C^*$, one has
$$\dim H^1(M,\LL_{\rho} \otimes f^*\LL_2)\ge -\chi(S)+1$$
for any local system $\LL_2 \in \T(S)$.

\end{cor}

 \begin{ex} \label{exfin1} ( The deleted $B_3$-arrangement )

We return to Example \ref{exMASTER} and apply the above discussion to this test case.
The corresponding mapping $f:M \to \C^*$ has $B=\{0, \infty \}$ and $C(f)=\{1\}$. Indeed, with obvious notation, we get the following divisors:
$D_0=L_1+L_4+2L_5$, $D_{\infty}=L_2+L_3+2L_7$ and $D_1=L_6+2L$ where $L:x+y-1=0$ is exactly the line from the $B_3$-arrangement that was
deleted in order to get Suciu's arrangement. Moreover, the associated fibration $f':M' \to S'$ in this case is just the fibration of the $B_3$-arrangement discussed in \cite{FY}, Example 4.6.

The line $L$ is the only new component that has to be deleted, therefore $m''(1)=2$.
Since none of the fibers $D_c$ (in our case there  is just one, for $c=1$) is multiple,  Theorem  \ref{thm5} implies that
$$T(f)=\Z/2\Z.$$
Let $\gamma_i =\gamma(L_i)$. We know that $\rho (\gamma_i)=\pm 1$ and to get the exact values we proceed as follows.
First note that we can choose  $\rho (\gamma_1)= 1$, since the associated torus is
$$f^*(\T(\C^*))=\{(t,t^{-1},t^{-1},t,t^2,1,t^{-2},1)~~ |~~ t \in \C^*\}.$$
(In fact the choice $\rho (\gamma_1)= -1$ produces the character $\rho'_W$ introduced in Example \ref{exMASTER}.)
Next let $\alpha= \sum _{i=1,7}\alpha_i \gamma_i \in H_1(M)$. Then $\alpha \in \ker f_*$ if and only if
\begin{equation} \label{e1}
\alpha_1+\alpha_4+2\alpha_5=\alpha_2+\alpha_3+2\alpha_7.
\end{equation}
In our case, the morphism $\theta:\ker f_* \to \Z/2\Z$ is given by $\alpha \mapsto \alpha_2+\alpha_3 -\alpha_6$
It follows that $\gamma_6 \in \ker f_*$ and $\theta(\gamma_6)=1 \in \Z/2\Z$. It follows that $\rho (\gamma_6)=-1$.

Next $\gamma_1+\gamma_2  \in \ker f_*$ and  $\theta(\gamma_1+\gamma_2    )=1 \in \Z/2\Z$. It follows that
$\rho(\gamma_1) \rho(\gamma_2)=-1$, i.e. $\rho(\gamma_2)=-1$. The reader can continue in this way and get 
the value of $\rho=\rho_W$ given above in  Example \ref{exMASTER}.

\end{ex}

 \begin{ex} \label{exfin2}    ( A more general example: the  $\A_m$-arrangement) 

Let $\A_m$ be the line arrangement in $\PP^2$ defined by the equation
$$x_1x_2(x_1^m-x_2^m)(x_1^m-x_3^m)(x_2^m-x_3^m)=0.$$
This arrangement is obtained by deleting the line $x_3=0$ from the complex reflection arrangement
associated to the full monomial group $G(3,1,m)$ and was studied in  \cite{C1} and in \cite{CDS}. The deleted $B_3$-arrangement studied above is obtained by taking $m=2$.

Consider the associated pencil
$$(P,Q)=(x_1^m(x_2^m-x_3^m),x_2^m(x_1^m-x_3^m)).$$
Then the set $B$ consists of two points, namely $(0:1)$ and $(1:0)$, and the set $C$ is the singleton
$(1:1)$, see for instance  \cite{FY}, Example 4.6. It follows that $m'(c)=1$, $m''(c)=m$ and hence
via Theorem \ref{thm5} we get
$$T(f)=\frac {\Z}{m\Z}.$$
Using Proposition \ref{prop6}, we expect $(m-1)$  1-dimensional components in $\V_1(M)$, and this is precisely what has been proved in  \cite{C1}, or in  Thm. 5.7 in \cite{CDS}.
There are $r=2+3m$ lines in the arrangement, and to describe these components we use the
coordinates 
$$(z_1,z_2,z_{12:1},...,z_{12:m}, z_{13:1},...,z_{13:m},  z_{23:1},...,z_{23:m})$$
on the torus $(\C^*)^r$ containing $\T(M)$. Here $z_j$ is associated to the line $x_j=0$, for $j=1,2$,
and $z_{ij:k}$ is associated to the line $x_i - w^kx_j$, where $i,j=1,3$, $k=1,...,m$, and
$w=\exp (2\pi {\sqrt -1}/m)$.
All the above 1-dimensional components have the same associated 1-dimensional subtorus
$$\T=f^*(\T(\C^*))=\{(u^m,u^{-m},1,...,1,u^{-1},...,u^{-1},u,...,u) ~~|~~ u \in \C^* \}$$
where $f:M \to \C^*$ is the morphism associated to the pencil $(P,Q)$, and each element
$1$, $u^{-1}$ and $u$ is repeated $m$ times.
The associated maximal isotropic subspace $E$ in $H^1(M,\C)$ is spanned by the 1-form
$$\omega= m \frac{dx_1}{x_1}-m \frac{dx_2}{x_2} -\sum_{k=1,m} \frac{dx_1-w^kdx_3}{x_1-w^kx_3  } + \sum_{k=1,m} \frac{dx_2-w^kdx_3}{x_2-w^kx_3  }.   $$
The patient reader may check that for any $\alpha \in H^1(M,\C)$, the vanishing
$\alpha \wedge \omega=0$ in $ H^2(M,\C)$ implies that $\alpha$ is a multiple of $\omega$
(this is the maximality condition in this 1-dimensional case).
Let $\gamma_c$ be an elementary loop about one line $L$  in the fiber $\CC_c$, with multiplicity 1, e.g. $L:x_1-x_2=0$. Similarly, let  $\gamma_b$ be an elementary loop about one line $L'$  in the fiber $\CC_b$, with multiplicity 1, where $b= \infty= (0:1)$, e.g. $L':x_2-x_3=0$.
And let  $\gamma_0$ be an elementary loop about one line $L_0$  in the fiber $\CC_0$, with multiplicity 1, where $0= (1:0)$, e.g. $L_0:x_1-x_3=0$.
One can show easily that

\noindent (i) the classes $[\gamma_c]$ and  $[\gamma_b + \gamma_0  ]$ in the group 
$T(f)$ are independent of the choices made;

\noindent (ii) $[\gamma_c] = - [\gamma_b + \gamma_0  ] $ is a generator of  $T(f).$

It follows that a torsion character $\rho \in \T(M)$ such that $\LL_{\rho}|F=\C_F$ and inducing a nontrivial character 
$\tilde \rho: T(f) \to \C^*$ is given by
$$\rho =(1,1, w^k,...,w^k, w^{-k}, ..., w^{-k}, 1,...,1)$$
for $k=1,...,m-1.$ Here $\tilde \rho ( [\gamma_c])=w^k$ and $\rho$ is normalized by setting the last
$m$ components equal to 1.

\end{ex}

 \begin{ex} \label{exfin3} (A non-linear arrangement)

Consider again the pencil $\CC: (P,Q)=(x_1^m(x_2^m-x_3^m),x_2^m(x_1^m-x_3^m))$ associated above to the $\A_m$-arrangement, for $m \ge 2$.
 We introduce the following new notation: $C=\{ (0:1),(1:0), (1:1)\}$. Let $B \subset \PP^1$ be a finite set such that
$|B|=k\ge 2$ and $B \cap C=\emptyset.$ Consider the curve arrangement in $\PP^2$ obtained by taking the union of the $3m$ lines given by
$$(x_1^m-x_2^m)(x_1^m-x_3^m)(x_2^m-x_3^m)=0$$
with the $k$ fibers $\CC_b$ for $b \in B$. Let $M$ be the corresponding complement and $f:M \to S:=\PP^1 \setminus B$
be the map induced by the pencil $\CC$. Then one has the following.

\noindent (i) $T(f)= \frac{\Z}{m\Z} \oplus \frac{\Z}{m\Z} \oplus   \frac{\Z}{m\Z}.$ Let $e_j$ for $j=1,2,3$ denote
the canonical basis of $T(f)$ as a $ \frac{\Z}{m\Z}$-module.

\noindent (ii) For a character $\tilde \rho: T(f) \to \C^*$, let $W_{\rho}=\LL_{\rho} \otimes f^*(\T(S))$ be the associated  component. Then $\dim W_{\rho}=k-1$ and for a local system $\LL \in W_{\rho}$ one has
$$\dim H^1(M,\LL) \ge k-2+\epsilon$$
where equality holds for all but finitely many  $\LL \in W_{\rho}$ and 
$$\epsilon = |\{j ~~|~~\tilde \rho (e_j) \ne 1\}|.$$
This shows that the various translates of the subtorus $\T_W= f^*(\T(S))$ have all the same dimension, but they are irreducible components of various characteristic varieties $\V_m(M)$, a fact not noticed before.

\end{ex}

\end{document}